\documentclass[12pt]{amsart}
\usepackage{lineno,amsmath,amsfonts}
\usepackage{amssymb,textcomp,epsfig,mathtools}
\usepackage{enumerate,caption}
\usepackage{xcolor}
\usepackage{url}
\usepackage{graphicx}
\usepackage{listings}
\makeatletter
\@namedef{subjclassname@2010}{%
 \textup{2010} Mathematics Subject Classification}
\makeatother

\newtheorem{conj}{Conjecture}

\theoremstyle{definition}

\numberwithin{equation}{section}
\def\NN{\mathbb{N}}

\usepackage{CJKutf8}
\newcommand{\hanyu}[1]{\begin{CJK}{UTF8}{gbsn}#1\end{CJK}}
\begin{document}

%%%%% To ease editing, for IMPAN journals add:

\baselineskip=17pt

%% In the running head, replace first names by initials and give an
%% abbreviation of the title.

\title[Zh\`\i-W\v{e}i S\={u}n's 1-3-5 Conjecture]{Report on\\[3mm]
  Zh\`\i-W\v{e}i S\={u}n's 1-3-5 Conjecture\\[3mm] and Some of Its
  Refinements}

\author[A. Machiavelo]{Ant\'onio Machiavelo} \address{Centro de
  Matem\'atica da Universidade do Porto and Mathematics Department of
  Faculdade de Ci\^encias do Porto, Rua do Campo Alegre 687, 4169-007
  Porto, Portugal} \email{ajmachia@fc.up.pt}

\author[R. Reis]{Rog\'erio Reis} \address{Centro de Matem\'atica da
  Universidade do Porto and Computer Science Department of Faculdade
  de Ci\^encias do Porto, Rua do Campo Alegre 1021/1055, 4169-007
  Porto, Portugal} \email{rvreis@fc.up.pt}

\author[N. Tsopanidis]{Nikolaos Tsopanidis} \address{Centro de
  Matem\'atica da Universidade do Porto, Rua do Campo Alegre 687,
  4169-007 Porto, Portugal} \email{tsopanidisnikos@gmail.com}

\date{15/05/2020}
%%%%%%%%%%%%%%%%%%%%%%%%%%%%%%%%%%%%%%%%%%%%%%%%%%%%%%%%%%%%%%%%%%%%%%%%

\begin{abstract}
  We report here on the computational verification of a refinement of
  Zh\`\i-W\v{e}i S\={u}n's ``1-3-5 conjecture'' for all natural
  numbers up to $105\,103\,560\,126$. This, together with a result of
  two of the authors, completes the proof of that conjecture.
\end{abstract}



\maketitle

%%%%%%%%%%%%%%%%%%%%%%%%%%%%%%%%%%%%%%%%%%%%%%%%%%%%%%%
\section{Introduction}

In a paper on refinements of Lagrange's four squares theorem,
Zh\`\i-W\v{e}i S\={u}n (\hanyu{孙智伟}) made the conjecture that any
$m\in\NN$ can be written as a sum of four squares, $x^2+y^2+z^2+t^2$,
with $x, y, z, t\in\NN_0$, in such a way that $x+3y+5z$ is a perfect
square. This is Conjecture 4.3(i) in \cite{sun2017refining}, and
Zh\`\i-W\v{e}i S\={u}n called it the ``1-3-5
conjecture''. Q\`{\i}ng-H\v{u} H\'ou (\hanyu{侯庆虎}) verified it up
to $10^{10}$: see \url{https://oeis.org/A271518}. We report here on
our computational verification of this conjecture for all
$m\leq 105\,103\,560\,126$.

This last number arose from the work \cite{MT2020}, in which it was
proved that the 1-3-5 conjecture is true for all numbers
$$m>\left(\frac{10}{\sqrt[4]{35}-\sqrt[4]{34}}\right)^4\simeq
105103560126.80255537.$$ Hence, the computation we are here reporting
on, together with the main result of \cite{MT2020}, completes the
proof that the 1-3-5 conjecture holds for all natural numbers.

Along the way, we have used another conjecture of S\={u}n, at his own
suggestion, Conjecture 4.9(ii) of \cite{sun2019restricted}, the part
whose content is as follows.
\begin{conj}[Zh\`\i-W\v{e}i S\={u}n]
  \label{conj:Sun19}
  Any positive integer can be written as $x^2+y^2+z^2+t^2$ with
  $x,y,z,t\in\NN_0$ such that $x+3y+5z$ is a square, and either $x$ is
  three times a square, or $y$ is a square, or $z$ is a square.
\end{conj}

As we will see below, we checked this conjecture for all natural
numbers up to $105\,103\,560\,126$, which implies the 1-3-5 conjecture
for the same range. We have noticed that, actually, for most numbers,
one may drop the possibility that $z$ is a square. Moreover, it now
seems that after some point on, all numbers have a representation as
in the statement of Conjecture~\ref{conj:Sun19}, but with
$x\in\{0,3\}$ or $y\in\{0,1\}$. This is the content of
Conjecture~\ref{newC} below.

%%%%%%%%%%%%%%%%%%%%%%%%%%%%%%%%%%%%%%%%%%%%%%%%%%%%%%% 
\section{The Modus Operandi}

We will say that a quadruple $(x,y,z,t)\in \NN_0^4$ is a \emph{1-3-5
  representation} of $m$ if $x^2+y^2+z^2+t^2 = m$ and $x+3y+5z$ is a
perfect square. Since it is clear that if a number $m$ has the 1-3-5
representation $(x,y,z,t)$, then $(4x,4y,4z,4t)$ is a 1-3-5
representation of $16m$, in order to verify the 1-3-5 conjecture, we
may disregard multiples of $16$.

Early on, during the first computations we made, it was found out
that, apparently, only the following $15$ numbers:
$$31, 43, 111, 151, 168, 200, 248, 263, 319, 456, 479, 871, 1752,
1864, 3544,$$
and their multiples by powers of 16, do not have a 1-3-5
representation $(x,y,z,t)$ where either $x$ is three times a square or
$y$ a square. We used this, together with Conjecture~\ref{conj:Sun19},
to speed up the search.

Furthermore, while testing the program's speed, and while running it
for values up to $10^7$, and then up $10^8$, it was noted that only
$123$ numbers have a special 1-3-5 representation requiring $x$ and
$y$ bigger than $4$, and that the last one of these was $779832$.
% 31,43,111,151,168,188,191,200,248,263,319,456,479,504,871,956,
% 1039,1199,1303,1399,1519,1528,1543,1663,1752,1864,2191,2463,2551,
% 2719,2872,3544,4039,4600,4639,5039,5128,5399,5416,5623,5719,5752,
% 6519,7484,8527,8824,9151,9431,9887,10136,10551,10991,11768,11848,
% 14536,14719,16279,16888,17560,17999,18424,18680,19576,19912,20984,
% 21224,22519,23551,24719,26088,26104,26927,30040,31992,32287,33016,
% 33208,33871,35704,35848,38872,39679,41512,49816,52728,61679,62248,
% 62712,69112,87544,87639,100344,102568,108024,108943,113608,119288,
% 125608,138744,146584,154552,156616,193912,202744,209959,239623,
% 323656,346312,366072,431608,432712,609784,767992,959992,1033720,
% 1211399,1393144,2035192,2161144,2272504,4833784,6558712,7779832
This observation motivated the introduction of a ``tolerance'' input
on the program, make it to exit whenever the list of numbers to be
checked was smaller than a certain size, and returning the list of
those numbers, which can then be checked individually in a faster way.

To tackle the verification of the 1-3-5 conjecture up to the required
number, $105\,103\,560\,126$, a program in the C programming language
was written that takes into account the above remarks, and runs as
follows. Firstly, it allocates the necessary memory for the range one
is checking, ignoring the multiples of $16$. Then, it looks for all
triplets $(x,y,z)$ such that $x+3y+5z$ is a square, and either $x$ is
three times a square or $y$ is a square. Each time such a triplet is
found, one removes from the appropriate memory the numbers
$m=x^2+y^2+z^2+s$, for all squares $s$ with $m$ in the desired
range. The program ends when the list of the remaining numbers has
size less than a prescribed number, which is part of the input.

Since what is wanted is to check, for every integer below a given
bound, if there is an additive decomposition in four squares, the
naive implementation of such check would have a $\mathcal{O}(n^3)$
complexity\footnote{All complexity considerations made here suppose
  that the cost of arithmetic operations for integers in the range
  considered has complexity $\mathcal{O}(1)$.}. But the simple
observation that we can check the existence of the fourth square
summand by subtracting the summation of the first three to the integer
that is tested, and check that this result is a square, lowers this
complexity to $\mathcal{O}(n^{\frac52})$. This is still a complexity
that makes the algorithm intractable for the desired bound. Since we
could not find a canonical ordering for the possible summands that
would allow to efficiently prune the search tree, the solution relied
in the classic space/time tradeoff. Thus, we represented the whole
integer search space as a bitmap, and with a
$\mathcal{O}(n^{\frac32})$ search could sweep this space and verify
the conjecture. As a matter of fact, and because, as already
mentioned, we can exhaust the whole search space with just a few
instances of the variable of the outer cycle, in practice the
algorithm finishes in a $\mathcal{O}(n)$ time.

The problem with this approach is that, because the size of the search
space is quite considerable, the memory space necessary to store the
bitmap representing the referred set was larger than the one available
in our laptops. Thus, the program does not try to cover the whole
space of considered integers in just one run, but splits this space in
various slices, that are searched independently. This has the
advantage of a ``parallelism for poor people'', running the program on
different slices in different laptops, but has the drawback that the
program for the higher slices, by the nature of the additive
decomposition, needs the same time to conclude as the same program
would need to run on an unsplitted space of integers. With laptops of
$16$GB of RAM, we splitted the search space in $11$ slices: the $i$-th
slice covering the range $[\,i\times 10^{10}, (i+1)\times 10^{10}\,]$,
for $0\leq i\leq 9$, and the $11$-th slice covering the remaing
numbers up to $105\,103\,560\,126$. In the process, the range
previously checked by Q\`{\i}ng-H\v{u} H\'ou was rechecked.

The code of the C program we used is given in Appendix~\ref{appC}. The
input consists of three numbers: the range over which one is checking;
the ``tolerance'', which is the size of the list of the numbers that
were not checked yet; and, finally, the interval one is checking.

%%%%%%%%%%%%%%%%%%%%%%%%%%%%%%%%%%%%%%%%%%%%%%%%%%%%%%% 
%%%%%%%%%%%%%%%%%%%%%%%%%%%%%%%%%%%%%%%%%%%%%%%%%%%%%%% 
\section{The results}

The different slices were distributed by several machines, and each
slice took between 2 to 3.5 days (depending on the machine, and on the
extra use that its owner was making of it). While working on slice 1,
it was noticed that only four numbers required the outter ``for''
cycle to go beyond $1$, and checking these four numbers was taking a
huge amount of time. Thus, the program was interruped, and restarted
with a tolerance of $10$, and here is the output:

\medskip

\small
\begin{verbatim}
  135Siever version 1.3
  Sieving 135 for 105103560126, with tolerance 10, 
  in the interval [10000000000,20000000000]

  0 9375000000
  1 1562500004
  4 4

  Done!! Lasting numbers: 4
  10234584952,11035927288,11051651704,14485001848
\end{verbatim}
\normalsize

\medskip

This output means that there are $9\,375\,000\,000$ numbers to be
checked (recall that one is ignoring the multiples of $16$), that
after the first run of the outter cycle (which looks for 1-3-5
representations $(x,y,z,t)$ where $x=3k^2$ or $y=k^2$, with
$k=0, 1, \ldots$), there remained only $1\,562\,500\,004$ numbers, and
that after the second run ($k=1$) only $4$ numbers are left. The
program then stops, and outputs those numbers.

These four $11$-digits-long positive integers were then checked using
the PARI/GP functions presented in appendix~\ref{appGP}, which uses an
algorithm to write a prime congruent to one modulo 4 as a sum of two
squares that is described by John Brillhart in \cite{Brilhart}. Using
those functions, one very quickly gets, for example (it is a random
algorithm), the following 1-3-5 representations:
\small
\begin{eqnarray*}
  (8524,9502,33094,94744) &\text{for}& 10234584952,\\
  (13438,32472,12774,98172)&\text{for}& 11035927288,\\
  (84720,34818,28982,42684)&\text{for}& 11051651704,\\
  (32742,93858,36824,56988)&\text{for}& 14485001848.
\end{eqnarray*}
\normalsize

\newpage

\noindent As a further example, we give here the output of the last
slice:

\bigskip

\small
\begin{verbatim}
  135Siever version 1.3
  Sieving 135 for 105103560126, with tolerance 10, 
  in the interval [100000000000,105103560126]

  0	4784587619
  1	797431269
  4	0

  Done!! Lasting numbers: 0

  132334.83 user 0.20 system 36:45:2 6elapsed 100%CPU 
  (0avgtext+0avgdata 623532maxresident)k 
  0inputs+8outputs (0major+155918minor)pagefaults 0swaps
\end{verbatim}
\normalsize

%%%%%%%%%%%%%%%%%%%%%%%%%%%%%%%%%%%%%%%%%%%%%%%%%%%%%%% 
\section{A new conjecture}

As a consequence of the computational results displayed above, we now
make the following conjecture.
\begin{conj}
  \label{newC}
  Any $m\in\NN$, that is not a multiple of $16$, with the exception
  of\, $31$, $43$, $111$, $151$, $168$, $200$, $248$, $263$, $319$,
  $456$, $479$, $871$, $1752$, $1864$, $3544$, can be represented as a
  sum of four squares, $x^2+y^2+z^2+t^2$, with $x,y,z,t\in\NN_0$ such
  that $x+3y+5z$ is a square, and either $x$ is three times a square,
  or $y$ is a square. Moreover, for $m > 14\,485\,001\,848$, one has a
  representation with $x\in\{0,3\}$ or $y\in\{0,1\}$, and
  (disregarding multiples of 16) exactly $\frac56$ of the numbers have
  a representation with $x=0$ or $y=0$, while the remainder $\frac16$
  have a representation with $x=3$ or $y=1$.
\end{conj}

\vspace{7mm}
%%%%%%%%%%%%%%%%%%%%%%%%%%%%%%%%%%%%%%%%%%%%%%%%%%%%%%%
\hrule
\appendix
\section{The C program}\label{appC}

\lstset{tabsize=2}

\tiny
\begin{lstlisting}
#include <stdio.h>
#include <stdlib.h>
#include <math.h>

#define VERSION "1.4"

#define MAX  10000000L
#define LIM  0
#define MAXS 10000L
#define MINS 0L
#define FULLSET (B64)0xfffefffefffefffe

typedef unsigned long B64;
typedef unsigned long Long;

typedef struct {
    B64 *map;
    Long min, max, nelements;
} BitMap;

BitMap map, squares;
B64 *masks;
Long max = MAX, mins=MINS, maxs=MAX;

B64* buildMasks(void){
    B64 *masks, *pt, val=(B64)1;
    int i;
    
    masks = (B64*)malloc(64*sizeof(B64));
    pt = masks;
    for(i=0; i<64; i++){
        *(pt++) = val;
        val = val << 1;
    }
    return masks;
}

void newBMapFull(BitMap *bmap, Long min, Long max){
    B64 *pt;
    Long nbytes, i, size;
    
    size = max - min + 1;
    nbytes = (size /(sizeof(B64)*8)) +1;
    bmap->map = (B64 *)malloc(sizeof(B64)*nbytes);
    bmap->max = max;
    bmap->min = min;
    bmap->nelements = max-(max/16) - (min-(min/16));
    pt = bmap->map;
    for(i=0;i<nbytes;i++) (*pt++) = FULLSET;
}

int memberP(BitMap *bmap, Long n){
    B64 byte, foo;
    Long rnumber;
    
    if (n > bmap->max || n < bmap->min) return 0;
    rnumber = n - bmap->min;
    byte =  rnumber/64;
    foo = rnumber - byte * 64;
    if(*(masks+foo) & *(bmap->map+byte)) return 1;
    return 0;
}

void removeM(BitMap *bmap, Long n){
    Long byte, foo, rnumber;

    if (n > bmap->max || n < bmap->min) return;
    rnumber = n - bmap->min;
    byte = rnumber/64;
    foo = rnumber - byte * 64;
    if(*(masks+foo) & *(bmap->map+byte)){
        (bmap->nelements)--;
        *(bmap->map + byte) = *(bmap->map+byte) & ~*(masks+foo);
    }
}

void addM(BitMap *bmap, Long n){
    B64 byte, foo;
    Long rnumber;

    if (n > bmap->max || n < bmap->min) return;
    rnumber = n - bmap->min;
    byte = rnumber/(64);
    foo = rnumber - byte * 64;
    if(!(*(masks+foo) & *(bmap->map+byte))){
        (bmap->nelements)++;
        *(bmap->map + byte) = *(bmap->map+byte) | *(masks+foo);
    }
}

void printM(BitMap *bmap){
    Long i, size;
    int j;
    
    size = bmap->max - bmap->min +1;
    for(i=0; i <= (bmap->max)/64;i++){
        if(*(bmap->map+i)){
            for(j=0; j<64; j++){
                if(i*64+j+(bmap->min) > bmap->max){
                    printf("\n");
                    return;
                }
                if(*(bmap->map+i) & *(masks+j)) printf("%lu ",i*64+j+bmap->min);
            }
        }
    }
}

void saveM(BitMap *bmap, int ord){
    FILE *file;
    char fname[100];
    Long i;
    int j;
    
    sprintf(fname, "c135-%d.csv",ord);
    file = fopen(fname,"w");
    for(i=0; i<= (bmap->max)/64;i++){
        if(*(bmap->map+i)){
            for(j=0; j<64; j++){
                if(i*64+j > bmap->max){
                    printf("\n");
                    return;
                }
                if(*(bmap->map+i) & *(masks+j)) fprintf(file,"%lu, ",i*64+j+bmap->min);
            }
        }
    }
    fclose(file);
}

int squarep(Long n){
    Long i;
    i = (int)(sqrt(n)+0.5);
    return i*i == n;
}

void dealWTriple(BitMap* map, Long i, Long j,Long k){
    Long foo, n=0, n2=0;
    if(squarep(i+3*j+5*k)){
        foo = i*i+j*j+k*k;
        while(1){
            if(foo + n2 > max) break;
            removeM(map, foo+n2);
            n2 += 2*(n++)+1;
        }
    }
}

int main(int argc, const char * argv[]) {
    Long i2=0, i=0, j, k, i4=0, lim=LIM;
    printf("135Siever version %s\n",VERSION);
    
    masks = buildMasks();
    if(argc == 4 ){
        lim = atol(argv[3]);
        mins =atol(argv[1]);
        maxs = atol(argv[2]);
        max = maxs;
        printf("Sieving 135, with tolerance %lu, in the interval [%lu,%lu]",lim,mins,maxs);
    } else {
        printf("Usage c135 lim min max\n");
        exit(-1);
    }
    newBMapFull(&map, mins, maxs);
    while(i4 <= max){
        printf("%lu\t%lu\n",i2,map.nelements);
        if(map.nelements <= lim){
            printf("\n");
            printf("Done!! Lasting numbers: %lu\n",map.nelements);
            printM(&map);
            printf("\n");
            exit(0);
        }
        for(j=0; j*j<= maxs-i4;j++){
            for(k=j; k*k<= (maxs-i4-j*j); k++){
                if(k*k+j*j+i2*i2 > maxs) break;
                dealWTriple(&map, 3*i2,j,k);
                dealWTriple(&map, 3*i2,k,j);
                dealWTriple(&map, j,i2,k);
                dealWTriple(&map, k,i2,j);
            }
        }
        i4 += 2*i2+1;
        i2 += 2*(i++)+1;
    }
    printf("\n");
    printM(&map);
    printf("\n Done! Lasting numbers: %lu\n",map.nelements);
    return 0;
}
\end{lstlisting}
\normalsize
%%%%%%%%%%%%%%%%%%%%%%%%%%%%%%%%%%%%%%%%%%%%%%%%%%%%%%%
\hrule
\bigskip

\section{The PARI/GP functions}\label{appGP}

\tiny
\begin{lstlisting}
/* Representation of a quaternion as a 4 x 4 matrix */

quat(a,b,c,d)=[a,b,c,d;-b,a,-d,c;-c,d,a,-b;-d,-c,b,a];

/* The Hurwitz units */

unid=[quat(1,0,0,0), quat(-1,0,0,0), quat(0,1,0,0), quat(0,-1,0,0), 
quat(0,0,1,0), quat(0,0,-1,0), quat(0,0,0,1), quat(0,0,0,-1), 
quat(1/2,1/2,1/2,1/2), quat(-1/2,-1/2,-1/2,-1/2), 
quat(1/2,1/2,1/2,-1/2), quat(-1/2,-1/2,-1/2,1/2),
quat(1/2,1/2,-1/2,1/2), quat(-1/2,-1/2,1/2,-1/2),
quat(1/2,-1/2,1/2,1/2), quat(-1/2,1/2,-1/2,-1/2),
quat(1/2,1/2,-1/2,-1/2), quat(-1/2,-1/2,1/2,1/2),
quat(1/2,-1/2,1/2,-1/2), quat(-1/2,1/2,-1/2,1/2), 
quat(1/2,-1/2,-1/2,1/2), quat(-1/2,1/2,1/2,-1/2),
quat(1/2,-1/2,-1/2,-1/2), quat(-1/2,1/2,1/2,1/2)];

/* Fast modular exponentiation */

expmod(a,e,m)={
  local(x,y,s,d);
  x=a; y=1; s=e;
  while(s,d=s%2;
    s=(s-d)/2; if(d,y=(y*x)%m); x=(x*x)%m);
  return(y);
  }

/* Imodp computes the solution of x^2 = -1 (mod p) with 0<x<p/2, for p=1
(mod 4) */

Imodp(p)={
  local(g,x);
  if(p%4<>1,return("not a valid prime!"));
  while(1,g=random(p);
    if(expmod(g,(p-1)/2,p)==p-1,x=expmod(g,(p-1)/4,p);
      if(x>p/2,x=p-x);return(x)));
  }

/* Euclid algorithm to compute gcd(a,b) but stopping at the first
remainder that is < sqrt(a) */

EuclSp(a,b)={
  local(r,x,z);
  r=a%b;
  if(r==1,return([b,1]));
  x=b;
  while(r>sqrt(a),z=x%r;x=r;r=z);
  return([r,x%r]);
  }

/* Writes p == 1 (mod 4) as a sum of two squares, using the algorithm
described in [1] */

PrimeSS(p)=EuclSp(p,Imodp(p));

/* A method to decompose an odd number as a sum of of either four
integer squares, or half integer squares, in a random way */

sum4sqF(a)={ 
  local(b,c,sq,x,y,z,tt,uu,vv,ct,s);
  if(a==1,return(quat(1,0,0,0)));
  while(1,
    b=random(2*sqrtint(a)-1)+1; if(b%2==0,b=b-1);
    c=4*a-b^2;
    sq=floor((sqrtint(c)+1)/2);
    j=random(sq-1)+1;
    z=(c-(2*j-1)^2)/2;
    if(isprime(z), v=PrimeSS(z); x=v[1];y=v[2];
      uu=x+y;vv=x-y;tt=2*j-1;
        s=matsolve([1,1,1,1;1,1,-1,-1;1,-1,1,-1;1,-1,-1,1],[b;uu;vv;tt]);
        return(quat(s[1,1],s[2,1],s[3,1],s[4,1]))
      ));
  }

/* Expanding the sum4sqF function to all natural numbers with results
only in the integers */

v2(n)={
  local(oddp);
  e=0;oddp=n;
  while(oddp%2==0,oddp=oddp/2;e=e+1);
  return(e);
  };

sum4sqFall(a)={
  local(r,e);
  e=v2(a);
  if(e==0,r=sum4sqF(a),a=a/2^e;r=(quat(1,1,0,0)^e)*sum4sqF(a));
  while(floor(r[1,1])<>r[1,1], 
         r=r*unid[random(length(unid))+1];
        );
  return(r);
  }

/* All possible permutations of the numbers 0,1,3,5 as 4D-vectors */

perm=List([[0, 1, 3, 5], [0, 1, 5, 3], [0, 3, 1, 5], [0, 3, 5, 1], 
[0, 5, 1, 3], [0, 5, 3, 1], [1, 0, 3, 5], [1, 0, 5, 3], [1, 3, 0, 5], 
[1, 3, 5, 0], [1, 5, 0, 3], [1, 5, 3, 0], [3, 0, 1, 5], [3, 0, 5, 1], 
[3, 1, 0, 5], [3, 1, 5, 0], [3, 5, 0, 1], [3, 5, 1, 0], [5, 0, 1, 3],
[5, 0, 3, 1], [5, 1, 0, 3], [5, 1, 3, 0], [5, 3, 0, 1], [5, 3, 1, 0]]);

/* rep135 gives a solution of the system 1-3-5, returning the
Lipschitz integer whose norm is the imput number, and the permutation
of 0135 with whom its inner product is a square */

rep135(a)={
  local(s,t,f,c);
  if(issquare(a), return([[sqrtint(a),0,0,0],[1,3,5,0]]));
  while(1, s=sum4sqFall(a);
    t=[abs(s[1,1]),abs(s[1,2]),abs(s[1,3]),abs(s[1,4])];
  for(i=1,length(perm),f=perm[i]*t~;
    if(issquare(f), return([t,perm[i]]))));
  }
\end{lstlisting}
\hrule
\bigskip

\normalsize
%%%%%%%%%%%%%%%%%%%%%%%%%%%%%%%%%%%%%%%%%%%%%%%%%%%%%%%
\section*{Acknowledgments}
  
The authors would like to thank Professor Zh\`\i-W\v{e}i S\={u}n
(\hanyu{孙智伟}) and Professor B\={o} H\'e (\hanyu{何波}) for their
kind feedback, helpful comments, and suggestions, and Gra\c{c}a Brites
and Vasco Machiavelo for letting their personal computers be used for
slices 3 and 4, respectively.

The authors would also like to acknowledge the financial support by
FCT --- Funda\c{c}\~ao para a Ci\^encia e a Tecnologia, I.P.---,
through the grants for Nikolaos Tsopanidis with references:\\
\centerline{PD/BI/143152/2019, PD/BI/135365/2017, PD/BI/113680/2015,}\\
and by CMUP --- Centro de Matem\'atica da Universidade do Porto ---,
which is financed by national funds through FCT under the project with
reference UID/MAT/00144/2020.
%%%%%%%%%%%%%%%%%%%%%%%%%%%%%%%%%%%%%%%%%%%%%%%%%%%%%%%
\bibliographystyle{alpha}
\bibliography{1-3-5}

%%%%%%%%%%%%%%%%%%%%%%%%%%%%%%%%%%%%%%%%%%%%%%%%%%%%%%%
\end{document}